\documentclass[12pt]{article}
\usepackage{amsmath,amssymb}
\usepackage{hyperref}
\newcommand{\punkt}[1]{\medskip\centerline{\large\bf #1}\medskip\par\noindent}
\newcommand{\RP}{{\bf R}^+}
\newcommand{\Th}{\par\medskip\noindent\large\bf}
\newcommand{\Rem}{\par\medskip\noindent\large}
\newcommand{\Proof}{\par\medskip\noindent\large}
\newcommand{\cR}{{\cal I}}
\newcommand{\N}{\bf\rm N}
\newcommand{\eps}{\varepsilon}
\begin{document}
\author{ Vyacheslav~S.~Rychkov \\\it Moscow Institute of Physics and Technology}
\title{Splitting of Volterra Integral Operators with Degenerate Kernels\noindent \footnote{This is author's translation of the Russian original published in {\it Investigations in the theory of differentiable functions of many variables and its applications. Part 17, Collection of articles}, Trudy Mat. Inst. Steklova, 214, Nauka, Moscow, 1997, 267–285; editorial translation published in: Proc. Steklov Inst. Math., 214 (1996), 260–278] \url{http://mi.mathnet.ru/eng/tm1040}.}
	\footnote{Work supported by the Russian Fund for Basic Research (project 96-01-00243).}
}
\date{}

\maketitle
\begin{abstract}
	Volterra integral operators with non-sign-definite degenerate kernels $A(x,t)= \sum_{k=0}^n A_k(x,t)$, $A_k(x,t)= a_k (x) t^k$, are studied acting from one weighted $L_2$ space on $(0,+\infty)$ to another. Imposing an integral doubling condition on one of the weights, it is shown that the operator with the kernel $A(x,t)$ is bounded if and only $n+1$ operators with kernels  $A_k(x,t)$ are all bounded. We apply this result to describe spaces of pointwise multipliers in weighted Sobolev spaces on $(0,+\infty)$.
\end{abstract}

\punkt{Introduction}
There exists a problem of studying weighted estimates of the form
$$
\int_0^\infty\left|v(x)\int_0^x A(x,t)f(t)\,dt\right|^p \,dx
\le C
\int_0^\infty|u(x)f(x)|^p\,dx
\eqno (1)
$$
from the viewpoint of finding necessary and sufficient conditions on 
$u(x)$ and $v(x)$,
under which inequality (1) holds for all functions $f$
with the finite r.h.s. of (1)
and a constant $C$ independent of $f$.
In the going back to Hardy [1,2]
case $A(x,t)=1$ the corresponding criterion was found in [3--5].

The first significant progress for $A(x,t)\ne1$ was achieved in [6--9], which investigated the case
$A(x,t)=(x-t)^\alpha$, $\alpha>0$. Later a number of works gradually extended the class of kernels
$A(x,t)$, for which it was possible to characterize $u(x)$ and $v(x)$ in (1).
At present the most general results were obtained, it seems, in~\cite{Oin}. This work finds a condition on $u(x)$ and $v(x)$, necessary and sufficient for
the validity of inequality (1) under the assumption that the kernel $A(x,t)$
is nonnegative and satisfy an additional condition of the form
$A(x,t)\asymp A(x,y)+A(y,t)$, $y\in(t,x)$,
allowing to compare kernel values in different points.

Section 1 of our work points out a new class of kernels, for which it is possible to characterize the weights in (1) for $p=2$ and an additional condition on $u(x)$.
The kernels of our class do not satisfy the above-mentioned conditions of applicability of known results.

In section 2 the results of section 1 are applied to the problem of describing the set of pointwise multipliers in some weighted Sobolev spaces. Sections 3--5 collect auxiliary results and proofs.

The author is deeply grateful to O.V. Besov and G.A. Kalyabin for valuable remarks, discussions and support.

{\bf Note added (June 2020)} See \cite{RycDokl} for a short presentation of these results without proofs, and \cite{Rychkov3} for generalizations to $p\ne2$.

\punkt{1. Weighted estimates of integral operators}
For a Volterra integral operator~$\cal A$ with a degenerate kernels of the form
$$
  ({\cal A}f)(x)=\int_0^x \left[\sum_{k=0}^n a_k (x)t^k \right]f(t) \,dt,
  \quad n \in \mbox{\bf N},  \eqno (1.1)
$$
we will study the possibility of an estimate
$$
 \| v {\cal A}f \|_2 \le C \|  uf \|_2.  \eqno  (1.2)
$$
Here $\|\cdot\|_2=\|\cdot\|_{L_2(\RP)}$; $\RP =(0,\infty)$;
$u,v$ are nonnegative on~$\RP$ functions (weights); constant~$C>0$
does not depend on $f$.

Denote by $L_{2,u}$ the weighted space of functions 
$f$ on $\RP$ with the norm $\|uf\|_2$. That (1.2) holds now means that
${\cal A}:L_{2,u}\to L_{2,v}$. We represent~$\cal A$ as a sum
${\cal A}=\sum_{k=0}^n {\cal A}_k$, where
$$
({\cal A}_k f)(x)=a_k(x)\int_0^x t^kf(t)\,dt.
$$
We will say that for the operator $\cal A$
when acting from $L_{2,u}$ into $L_{2,v}$ the {\it splitting} takes place, if
$$
{\cal A}:L_{2,u}\to L_{2,v}\Longleftrightarrow{\cal A}_k:L_{2,u}\to L_{2,v},
\quad k=0\ldots n.
$$

By $B_\delta,\delta\ge0$, denote the set of positive locally summable
on~$\RP$ functions~$w$, satisfying with some constant~$D(w)$ the integral doubling condition
$$
\int_\Delta w(x) \,dx \le D(w) \int_{\frac{1}{2}\Delta} w(x)\,dx
$$
for any interval $\Delta\subset\RP$ of length $|\Delta|\ge\delta$,
where $\frac{1}{2}\Delta$ is a twice smaller interval with the same center.

The following theorem, given in two equivalent formulations, is the central result of our work.

{\Th Theorem 1.1.}
{\it Let $u^{-2}\in B_\delta$ for some $\delta\ge 0$.
If $\delta>0$, then assume in addition $a_k v \in L_2 (0,r)\,\forall r>0$, $k=0\ldots n-1$.
Then for the inequality {\rm (1.2)} to be satisfied it is necessary and sufficient that
$$
S_k=\sup_{r>0} \|a_k v\|_{L_2 (r,\infty)}\cdot
		   \|x^k u^{-1}\|_{L_2 (0,r)} <\infty,\quad  k=0\ldots n.
                   \eqno  (1.3)
$$

}
{\Th Theorem 1.1$'$.}
{\it Assume the conditions of Theorem {\rm 1.1}. Then for the operator $\cal A$
acting from $L_{2,u}$ into $L_{2,v}$ the splitting takes place.
}

{\Rem Remark 1.} Our method of proof of Theorem 1.1
gives the following estimate for the norm of ${\cal A}$ (or, which is the same, the smallest constant
$C$ in inequality (2)):
$$
c_1\sum_{k=0}^nS_k\le\|{\cal A}\|_{L_{2,u}\rightarrow L_{2,v}}\le
c_2\sum_{k=0}^nS_k.
$$
Constant $c_1$ here depends on $n$, $\delta$, $D(u^{-2})$,
as well as (if $\delta>0$) on the quantity
$\sum_{k=0}^{n-1} \|a_k v\|_{L_2 (0,r_0)}\cdot\|x^ku^{-1}\|_{L_2 (0,r_0)}$,
where $r_0$ is determined by $n$, $\delta$, $D(u^{-2})$.
Constant $c_2$ is universal.

{\Rem Remark 2.} Consider the adjoint operator ${\cal A}^*$:
$$
({\cal A}^*f)(x)=\int_x^\infty \left[\sum_{k=0}^n x^ka_k (t)\right]f(t) \,dt.
$$
Since $\left(L_{2,u}\right)^*=L_{2,u^{-1}}$, we have
$ {\cal A}:L_{2,u}\to L_{2,v}
\Longleftrightarrow {\cal A}^*:L_{2,v^{-1}}\to L_{2,u^{-1}},
$
and so, under the assumptions of Theorem 1.1, Eq. (1.3) is also necessary and sufficient for having the inequality
$ \|u^{-1}{\cal A}^*f\|_2\le C\|v^{-1}f\|_2.$
\medskip

Condition $u^{-2}\in B_\delta$ of Theorem 1.1 allows to include many interesting cases.
E.g. weight
$u(x)=(1+x)^\alpha\log^\beta(2+x)$ satisfies this condition with $\delta=0$
for any $\alpha,\beta\in\bf R$. Nevertheless, it is natural to ask to what extent this condition is essential for the validity of Theorem 1.1. The rest of this section is devoted to clarifying this question. The available results are closely related with the paper \cite{Kuf}.

Consider the Riemann-Liouville integral operator:
$$
(\cR^{(\alpha)}f)(x)=
\frac{1}{\Gamma(\alpha)}\int_0^x(x-t)^{\alpha-1}
f(t)\,dt,\quad \alpha\ge1.
$$
For $\alpha\in\hbox{\bf N}$, operator $\cR^{(\alpha)}$ is an operator of the form (1.1).
Let us focus on $\alpha=2$ and
represent $\cR^{(2)}$ as a sum of two operators:
$$
\cR^{(2)}=\cR^{(2)}_0+\cR^{(2)}_1,
$$
$$
\cR^{(2)}_0f(x)=x\int_0^xf(t)\,dt,\quad
\cR^{(2)}_1f(x)=-\int_0^xtf(t)\,dt.
$$

{\Th Assertion 1.2.}
{\it For operator $\cR^{(2)}$, acting in the space
$L_{2,e^{-x}}$, splitting does not take place. Namely,
$$
   \cR^{(2)}:L_{2,e^{-x}}\to L_{2,e^{-x}},
$$
while
$$
 \cR^{(2)}_i:L_{2,e^{-x}}\not\to L_{2,e^{-x}},\quad i=0,1.
$$

}
This result is basically a reformulation of example 1 in \cite{Kuf}. It shows that condition
$u^{-2}\in B_\delta$ of Theorem 1.1 is important (clearly,
$u(x)=e^{-x}$ does not satisfy this condition for any
$\delta\ge0$).

Nevertheless, condition
$u^{-2}\in B_\delta$ is not, generally speaking, necessary for splitting. E.g. for
operators $\cR^{(\alpha)}$ it can be replaced by a weaker condition:
$$
\int_0^{2r}u^{-2}\,dx\le D\int_0^r u^{-2}\,dx\quad\mbox{for all }
r\ge\delta\ge0.         \eqno (1.4)
$$
Namely, the following result holds (basically obtained in \cite{Kuf}, although not formulated there explicitly).

{\Th Assertion 1.3.}
{\it Let the weight $u$ satisfy with some constants
$D$, $\delta$ condition {\rm (1.4)}.
If $\delta>0$, then let in addition
$x^{\alpha-1}v\in L_2(0,r)\forall r>0$. Then for $\alpha\ge 1$ to have the inequality 
$$   \| v \cR^{(\alpha)}f \|_2 \le c \|  uf \|_2
\eqno (1.5)$$
it is necessary and sufficient that
$$
\sup_{r>0}\|x^{\alpha-1}v\|_{L_2(r,\infty)}\cdot
\|u^{-1}\|_{L_2(0,r)}<\infty. \eqno(1.6)
$$
For $\alpha\in \N$ this is equivalent to splitting for operator $\cR^{(\alpha)}$.
}
\medskip

\punkt{2. Pointwise multipliers in weighted Sobolev spaces}
Consider on $\RP$ the weighted Sobolev space
$W=W_{2,u}^{(l)}$ with the norm $\|f\|_W=\|f\|_{L_2(0,1)}+
\|f^{(l)}u\|_2$. For this norm, when the norm of function itself is taken only over an initial interval 
of $\RP$, all polynomials of degree $\le l-1$ belong to $W$. Spaces $W_{2,u}^{(l)}$
were introduced and studied in \cite{Kud}, which used an equivalent norm $\sum_{k=0}^{l-1}|f^{(k)}(0)|+\|f^{(l)}u\|_2$.

Function $\varphi$ is called a (pointwise) multiplier from
one Sobolev space ${}^1W$ to another ${}^2W$
if $\varphi f\in {}^2W\;\forall f\in{}^1W$. The space of 
such multipliers is denoted $M({}^1W\to {}^2W)$.

Various aspects of the theory of multipliers in spaces of differentiable functions were studied e.g. in the book~\cite{MSh}.

We are considering the problem of describing the space
$M(W_{2,u}^{(l)}\to W_{2,v}^{(m)})$, $m\le l$, denoted for brevity $M(u,l;v,m)$. In connection with this problem one should mention the work \cite{Kal}, which described multipliers in the Sobolev space on ${\bf R}^n$ with the norm
$\|f\|_{L_p(B(0,1))}+\sum_{|\alpha|=l}\|D^{(\alpha)}f\|_{L_p({\bf R}^n)}$
for the case $p>n$.

The first 3 assertions of this section are slight generalizations of author's results published in \cite{Ryc}.

{\Th Lemma 2.1.}
{\it Let function $g$ on $\RP$ be such that
$$
g^{(k)}(0)=0,\quad k=0\ldots l-1. \eqno (2.1)
$$
Then the following two equations hold:
$$
(\varphi g)^{(m)}(x)
=\frac1{(l-1)!}\sum_{k=0}^{l-1}C_{l-1}^k
\left(\varphi x^k\right)^{(m)}\int_0^x(-t)^{l-k-1}g^{(l)}(t)\,dt,
\quad m<l,
\eqno (2.2)
$$
$$
(\varphi g)^{(l)}(x)
=\varphi(x)g^{(l)}(x)+\frac1{(l-1)!}\sum_{k=0}^{l-1}C_{l-1}^k
\left(\varphi x^k\right)^{(l)}\int_0^x(-t)^{l-k-1}g^{(l)}(t)\,dt.
\eqno (2.3)
$$
}

\noindent The next lemma concerns the case $m=l$.
{\Th Lemma 2.2.}
{\it Let $u^{-1},v^{-1}\in L_2(0,r)\,\forall r>0$. Then
$\|\varphi v u^{-1}\|_{L_\infty(\RP)}<\infty$ for any function $\varphi\in M(u,l;v,l)$.}

\medskip\noindent From Lemmas 2.1 and 2.2 one easily derives

{\Th Theorem 2.3.} {\it Let $(1+x^{l-1})u^{-1}\in L_2(\RP)$,
$v^{-1}\in L_2(0,r)
\,\forall r>0$. Then the space
$M(u,l;v,m)$, $m\le l$, consists of those and only those $\varphi$ which satisfy the conditions
$$
\|(\varphi x^k)^{(m)}v\|_2<\infty,\quad k=0\ldots l-1,
\eqno (2.4)
$$
and in the case $m=l$ additionally
$$
\|\varphi v u^{-1}\|_{L_\infty(\RP)}<\infty.
\eqno (2.5)
$$

}

{\Rem Remark.} Theorem 2.3 states, roughly speaking, that
(in its conditions) to check whether $\varphi$ belongs to
$M(u,l;v,m)$ one should see how
multiplication by $\varphi$ acts on polynomials in $W_{2,u}^{(l)}$.
In other words, we have a weight effect: the growth of $u$
at $\infty$ implied by the condition $(1+x^{l-1})u^{-1}\in L_2$ leads to the fact that the functions 
of the space $W_{2,u}^{(l)}$ ``differ little'' from polynomials.
Note in this regard a result from \cite{Kud}: for
$(1+x^{l-1+\varepsilon})u^{-1}\in L_2$, $\eps>0$, for each function
$f\in W_{2,u}^{(l)}$
there exists a polynomial $P$ of degree $l-1$, to which it stabilizes in the sense that $\lim_{x\to\infty}(f(x)-P(x))^{(k)}=0$, $k=0\ldots l-1$.

The next result is deeper than Theorem 2.3 and needs Theorem 1.1 for its proof.

{\Th Theorem 2.4.}
{\it Let $u^{-2}\in B_\delta$, $v^{-1}\in L_2(0,r)
\,\forall r>0$. Then the space
$M(u,l;v,m)$, $m\le l$, consists of those and only those $\varphi$ which satisfy the conditions
$$
\|(\varphi x^k)^{(m)}v\|_2<\infty,\quad k=0\ldots l-1,
\eqno (2.6)
$$
$$
\begin{array}{c}
\mathop{\sup}\limits_{r>0}\|(\varphi x^k)^{(m)}v\|_{L_2(r,\infty)}\cdot
\|x^{l-k-1}u^{-1}\|_{L_2(0,r)}<\infty,\\[0.4mm]
k=0\ldots l-1,
\end{array}
\eqno  (2.7)
$$
and in the case $m=l$ additionally
$$
\|\varphi v u^{-1}\|_{L_\infty(\RP)}<\infty.
\eqno (2.8)
$$

}

Function $u(x)=(1+x)^\alpha$ for $\alpha<l-1/2$ is an example of a weight
allowed in Theorem 2.4, but not satisfying conditions of Theorem 2.3. The opposite example is provided by
$u(x)=e^x$: in this case it is Theorem 2.4 which is not applicable, and we must use
Theorem 2.3. Finally note that $u(x)=e^{-x}$ is not covered by any of these theorems; describing the corresponding multiplier spaces is a problem for the future.

\newpage
\punkt{3. Auxiliary results}
By $\Delta$, $\Delta_1$, $\Delta'\ldots$ we denote intervals in~$\RP$,
by $a\Delta$, $a>0$, the interval of length $a|\Delta|$ having the same center as $\Delta$. $\square$ denotes the end of proof.

{\Th Lemma 3.1.}{\it Let $w\in B_\delta$, and let function $\psi\ge 0$ satisfy the condition
$$
\sup_\Delta\psi\le c\inf_{\frac{1}{2}\Delta}\psi,\quad\mbox{if  }
|\Delta|\ge\delta.
$$
Then $\psi w\in B_\delta$. In particular, $x^\gamma w\in B_\delta\,
\forall\gamma>0$.

}
{\Proof Proof} is obvious. $\square$

The r.h.s. inequality of the next lemma is not surprising, while the l.h.s. one shows that the class $B_\delta$ is  more narrow than it could seem from the first glance.

{\Th Lemma 3.2.}{\it Let $w\in B_\delta$. Then there exist constants
$\alpha,\beta,A,B>0$ such that
$$
A\left(\frac{|\Delta_2|}{|\Delta_1|}\right)^\alpha\le
\frac{\int_{\Delta_2}w\,dx}{\int_{\Delta_1}w\,dx}\le
B\left(\frac{|\Delta_2|}{|\Delta_1|}\right)^\beta,\quad
\mbox{if  }
\Delta_1\subset\Delta_2, |\Delta_1|\ge\delta.
\eqno (3.1)
$$

}
{\Proof Proof.} Let us first prove the estimate
$$
\int_{2\Delta\cap\RP}w\,dx\le E\int_\Delta w\,dx,\quad
|\Delta|\ge\delta,
$$
with the constant $E=D(w)^2$. Obviously we only need to consider
$2\Delta\not\subset\RP$.
Let $\Delta=[a,a+2\varepsilon]$,
 $a<\varepsilon$. Then $2\Delta=[a-\varepsilon,
a+3\varepsilon]$, $2\Delta\cap\RP=[0,a+3\varepsilon]$.
Now
$$
\frac{1}{4}(2\Delta\cap\RP)=
\left[\frac{3}{8}(a+3\varepsilon),\frac{5}{8}(a+3\varepsilon)\right]
\subset\Delta,
$$
and therefore
$$
\int_{2\Delta\cap\RP}w\,dx\le
D(w)^2\int_{\frac{1}{4}(2\Delta\cap\RP)}w\,dx\le
D(w)^2\int_\Delta w\,dx.
$$
Now, to prove the r.h.s. inequality in (3.1) note that
$|2\Delta\cap\RP|\ge\frac{3}{2}|\Delta|$ for arbitrary $\Delta$.
Therefore if we choose an integer $N$ from the condition
$(3/2)^{N-1}<|\Delta_2|/|\Delta_1|\le(3/2)^N$, then after applying to
$\Delta_1$ $N$ consecutive operations $\Delta\mapsto 2\Delta\cap\RP$
we will get an interval covering $\Delta_2$. this implies that
$$
\frac{\int_{\Delta_2}w\,dx}{\int_{\Delta_1}w\,dx}\le E^N\le
E\left(|\Delta_2|/|\Delta_1|\right)^\beta
$$
for $\beta=\log_{3/2}E$, proving the r.h.s. inequality in (3.1).

To prove the l.h.s. inequality in (3.1), consider first the case when $\Delta_1$ and
$\Delta_2$ have the same left endpoint. Let $\Delta_1=[a,a+\varepsilon]$.
Consider the interval $\widetilde\Delta_1=[a+\varepsilon,a+2\varepsilon]$.
Then $\Delta_1\subset 4\widetilde\Delta_1$, therefore
$$
\int_{\Delta_1}w\,dx\le\int_{4\widetilde\Delta_1}w\,dx\le
E^2\int_{\widetilde\Delta_1}w\,dx=
E^2\left[\int_a^{a+2\varepsilon}-\int_a^{a+\varepsilon}w\,dx\right],
$$
and so
$$
\int_a^{a+2\varepsilon}w\,dx\ge\left(1+\frac{1}{E^2}\right)
\int_a^{a+\varepsilon}w\,dx.
$$
Applying this inequality $N$ times, where $2^N\le
|\Delta_2|/|\Delta_1|<2^{N+1}$,
we get
$$
{\int_{\Delta_2}w\,dx\over \int_{\Delta_1}w\,dx}\ge
\left(1+\frac{1}{E^2}\right)^N\ge
\left(1+\frac{1}{E^2}\right)^{-1}
\left(|\Delta_2|/|\Delta_1|\right)^\alpha,
$$
where $\alpha=\log_2\left(1+1/E^2\right)$.

For the general relative position of
$\Delta_1,\Delta_2$ consider intervals $\Delta_2'$ and $\Delta_2''$ such that
1)~$\Delta_2'\cup\Delta_2''=\Delta_2$, $\Delta_2'\cap\Delta_2''=\Delta_1$;
 2)~$\Delta_1$ and $\Delta_2'$ have the same right endpoint; 3) $\Delta_1$ and $\Delta_2''$ have the same left endpoint.
Then
$$
\left(\frac{|\Delta_2|}{|\Delta_1|}\right)^\alpha =
\left(\frac{|\Delta_2'|+|\Delta_2''|-|\Delta_1|}{|\Delta_1|}\right)^\alpha\le
2^\alpha\left(\left(\frac{|\Delta_2'|}{|\Delta_1|}\right)^\alpha+
	\left(\frac{|\Delta_2''|}{|\Delta_1|}\right)^\alpha\right)\le
$$
$$
\le
{2^\alpha}/A\left(\frac{\int_{\Delta_2'}w\,dx}{\int_{\Delta_1}w\,dx}+
					   \frac{\int_{\Delta_2''}w\,dx}{\int_{\Delta_1}w\,dx}\right)
={2^\alpha}/A\frac{\int_{\Delta_2'}w\,dx+\int_{\Delta_2''}w\,dx}
				   {\int_{\Delta_1}w\,dx}\le
{2^{\alpha+1}}/A\frac{\int_{\Delta_2}w\,dx}{\int_{\Delta_1}w\,dx}.
$$
$\square$

{\Th Lemma 3.3.} {\it Let ${P}_n(x)$ be an arbitrary degree $n$ polynomial, $\Delta_1\subset\Delta_2$ intervals in {\bf R}. Then for some constant $c=c(n)$
$$
\max_{\Delta_2}|{P}_n(x)|\le c
\left(\frac{|\Delta_2|}{|\Delta_1|}\right)^n
\max_{\Delta_1}|{P}_n(x)|.
$$

}
{\Proof Proof} of this undoubtedly known fact is given for completeness. Iterating Markov's inequality
$$
\max_{\Delta_1}|{P}_n'(x)|\le
\frac{2n}{|\Delta_1|}\max_{\Delta_1}|{P}_n(x)|,
$$
we get an estimate for the derivative of order $k$:
$$
\max_{\Delta_1}|{P}_n^{(k)}(x)|\le
\frac{c_1(n)}{(|\Delta_1|)^k}\max_{\Delta_1}|{P}_n(x)|,
\quad k\le n.
\eqno (3.2)
$$
Choose a point $x_0\in\Delta_1$ and write the Taylor expansion of
${P}_n(x)$ around $x_0$:
$$
{P}_n(x)=\sum_{k=0}^n\frac{(x-x_0)^k}{k!}{P}_n^{(k)}(x_0).
$$
Taking into account (3.2) this implies
$$
\max_{\Delta_2}|{P}_n(x)| \le
c_2(n)\sum_{k=0}^n\left(\frac{|\Delta_2|}{|\Delta_1|}\right)^k
\max_{\Delta_1}|{P}_n(x)|
$$
$$\le c(n)\left(\frac{|\Delta_2|}{|\Delta_1|}\right)^n
\max_{\Delta_1}|{P}_n(x)|.
\quad\square
$$

The next lemma is of a mostly technical character and describes some properties of polynomials orthogonal with weight $B_\delta$. Part (i)
asserts, roughly speaking, that roots of such polynomials cannot get too close to each other. 
Part (ii) is more particular and it will play an important role in the proof of Lemma 3.5.

{\Th Lemma 3.4.} {\it Let $v$, $w\in B_\delta$. Let ${P}_{n,r}$ be the $n$-th
polynomial of the orthogonal polynomial system with weight $v$ on the interval $[0,r]$, i.e.
$$
\int_0^rt^k{P}_{n,r}(t)v(t)\,dt=0,\quad k=0\ldots n-1. \eqno (3.3)
$$
{\rm (i)} Let $0<t_{r,1}<\ldots<t_{r,n}<r$ be the roots $P_{n,r}$, breaking up
$[0,r]$ into $(n+1)$ intervals $\Delta_{r,0},\Delta_{r,1},\ldots,\Delta_{r,n}$.
There exists a constant $\varepsilon=\eps(D(w),\delta,n)>0$ such that
for all $r\ge r_0=\delta/\eps$
$$
|\Delta_{r,j}|\ge\varepsilon r,\quad j=0\ldots n.
$$
{\rm (ii)}
Let us impose normalization $P_{n,r}(0)=1$. Consider the polynomial
$Q_{n,r}(t)=1-P_{n,r}(t)$. There exist
constants $\beta\in(0,1)$, $\gamma>1/\beta-1$, $r_0\ge0(=0$ for $\delta=0)$
such that for each $r\ge r_0$ one can select a set
$A_r\subset[0,r]$, for which
$$
\int_{A_r}w\,dt=\beta\int_0^r w\,dt,  \eqno (3.4)
$$
$$
\int_{[0,r]\backslash A_r}Q_{n,r}^2(t)w(t)\,dt\ge \gamma
\int_{A_r}Q_{n,r}^2(t)w(t)\,dt.     \eqno  (3.5)
$$
Constants $\gamma,\beta,r_0$ depend only on
$D(v), D(w),\delta, n$.
}

{\Proof Proof.} (i)
We will be omitting the lower index $r$ from the notation.
That the roots of 
${P}_n$ are simple, real and located on $(0,r)$ is not an additional requirement, but follows from (3.3), as is shown in the theory of orthogonal polynomials. Furthermore, for each
$r>0$ there exists, clearly, such a $j_0\in\{0,\dots,n\}$, that
$$
|\Delta_{j_0}|\ge\frac{r}{n+1}. \eqno (3.6)
$$
For $j\ne j_0$ let us consider the polynomial $R(t)=\prod (t-t_k)$, where the product is taken over $k\in\{1,...,n\}\backslash\{j,j+1\}$,
and let us use orthogonality of $R$ and ${P}_n$.
Let for definiteness $j_0\ne1$, $j=1$. Then
$$
\int_0^r(t-t_1)(t-t_2)\prod_{k=3}^n(t-t_k)^2w(t)\,dt=0,
$$
whence
$$
\int_{\Delta_1}(t-t_1)(t_2-t)\prod_{k=3}^n(t-t_k)^2w(t)\,dt\ge
$$
$$
\ge\int_{\Delta_{j_0}}(t-t_1)(t-t_2)\prod_{k=3}^n(t-t_k)^2w(t)\,dt.
$$
By (3.6), the polynomial under the last integral sign is not less than
$c(n)r^{2n-2}$ on $\frac{1}{2}\Delta_{j_0}$, while on  $\Delta_1$ it does not exceed
$r^{2n-2}$, therefore
$$
\int_{\Delta_1}w\,dt\ge c(n)\int_{\frac{1}{2}\Delta_{j_0}}w\,dt.
\eqno (3.7)
$$
Let $r_1=2(n+1)\delta$, then
$|\frac{1}{2}\Delta_{j_0}|\ge\delta$
for $r\ge r_1$, and we can continue inequality (3.7)
with the help of the r.h.s. inequality from (3.1) (Lemma 3.2):
$$
\int_{\Delta_1}w\,dt\ge c_0\int_0^rw\,dt.
\eqno (3.8)
$$
Let $r_0=\delta/\varepsilon\ge r_1$, where  $\varepsilon\in(0,(2n+2)^{-1})$
is taken so small, that for $r\ge r_0$ and $\Delta\subset[0,r]$,
$|\Delta|=\varepsilon r$, the l.h.s. inequality from (3.1) gives
$$
\int_0^rw\,dt\ge1/c_0\int_\Delta  w\,dt.
$$
From here and from (3.8) for $r\ge r_0$ we get
$$
\int_{\Delta_1}w\,dt\ge\int_\Delta w\,dt,\quad
\mbox{if }\Delta\subset[0,r],
|\Delta|=\varepsilon r.
$$
Therefore, $|\Delta_1|\ge\varepsilon r$, Q.E.D.

(ii)
On each interval $\Delta_j$ (see (i)) let us find the (clearly unique) point
$z_j$ where $|P_n(t)|$ attains the maximum on this interval. Clearly,
$z_0=0$, $z_n=r$, $P_n'(z_j)=Q_n'(z_j)=0$, $j=1\ldots n-1$.
Lemma 3.3 and (i) imply the inequality
$$
\max_{[0,r]}|P_{n}(t)|\le
c(n)\varepsilon^n\max_{\Delta_{j}}|P_{n}(t)|,
\quad r\ge r_0,
$$
from where
$$
|P_n(z_j)|=\max_{\Delta_{j}}|P_{n}(t)|\ge c_1,
\quad r\ge r_0.
\eqno (3.9)
$$

Denote
$F_j=[z_{j-1},z_j]$, $j=1\ldots n$, and consider the intervals
$$
I_j(\alpha)=\{t\in F_{j}:|Q_{n}(t)|\le
\min_{F_{j}}|Q_{n}(\cdot)|+\alpha\}.
$$
Choose $\alpha_j$ so that
$$
\int_{I_j(\alpha_j)}w\,dt=\beta\int_{F_j}w\,dt
$$
with some yet undetermined constant $\beta\in(0,1)$, and put
$$
A_r=\mathop{\cup}\limits_jI_{j}(\alpha_j).
\eqno (3.10)
$$
Equation (3.4) will be, clearly, satisfied. Part (ii) will be proved if
we can show that $\beta$ can be chosen so that the inequality
$$
\frac{\int_{F_j\backslash I_j(\alpha_j)}Q_{n}^2(t)w(t)\,dt}
{\int_{I_j(\alpha_j)}Q_{n}^2(t)w(t)\,dt}\ge \gamma
\eqno (3.11)
$$
holds for each $j$ with a constant $\gamma>1/\beta-1$.

Consider separately the two cases: 1) $\min_{F_j}|Q_n|=0$ and 2) $\min_{F_j}|Q_n|>0$.
In the first case (it is realized e.g. for $j=1$) for $\alpha_j<1$
we have
$$
\mbox{l.h.s. of (3.11) }\ge
{\int_{F_j\backslash I_j(1)}w\,dt \over
\alpha_j^2\int_{I_j(\alpha_j)}w\,dt}=
1/(\alpha_j^2\beta){\int_{F_j\backslash I_j(1)}w\,dt \over
\int_{F_j}w\,dt}.
$$
The closure of the set $F_j\backslash I_j(1)$ contains the interval $F_j^0=\{t\in F_j:P_n(t)\le0
\}$.
From $|\Delta_{j}|\ge\varepsilon r$ it is easy to deduce an analogous
inequality for the intervals $F_{j}^0$. Let us increase if necessary $r_0$
so that for $r\ge r_0$ we have $|F_j^0|\ge\delta$. Then by Lemma 3.2 we have
$\left.\int_{F_j\backslash I_j(1)}w\,dt\right/\int_{F_j}w\,dt\ge c$,
and inequality
(3.11) is satisfied with $\gamma=c/(\alpha_j^2\beta)$.
Note that $\alpha_j=O(|I_j(\alpha_j)|)=O(\beta)$.
[The first equation uses the following from Markov's inequality and Lemma 3.3 estimates $Q_n'(x)=O(1/r)$ on $\Delta_0$ $\Rightarrow O(1)$ on $[0,r]$;
the second equation follows from Lemma 3.2.] Now it's clear that for sufficiently 
small $\beta>0$ we will have
$\gamma>1/\beta-1$.

In the second case we argue as follows. Polynomial $Q_n$ varies monotonically and does not vanish on
$F_j$. For definiteness assume that it's positive and increasing.
From (3.9) it follows that $0<Q_n(z_{j-1})=\min_{F_j}Q_n\le1-c_1$. Furthermore we can find a $\xi\in(z_{j-1},z_j)$ such that
$$
\int_{F_j}Q_n^2w\,dt=Q_n^2(\xi)\int_{F_j}w\,dt,
$$
From the simple estimates
$$
\int_{F_j}Q_n^2w\,dt\ge(\min_{F_j}Q_n)^2\int_{F_j\backslash F_j^0}w\,dt+
\int_{F_j^0}w\,dt=
$$
$$
=(\min_{F_j}Q_n)^2\int_{F_j}w\,dt+
=1-(\min_{F_j}Q_n)^2)\int_{F_j^0}w\,dt\ge
$$
$$
\ge\left[(\min_{F_j}Q_n)^2+c(1-(1-c_1)^2)\right]\int_{F_j}w\,dt
$$
it follows that
$Q_n(\xi)\ge\min_{F_j}Q_n+\theta$, $\theta>0$.
Denote by $\eta$ the right endpoint of $I_j(\alpha_j)$, i.e.
$I_j(\alpha_j)=[z_{j-1},\eta]$.
Monotonicity of $Q_n$ and the definition of $\xi$ imply the inequality
$$
\int_\eta^{z_j}Q_n^2w\,dt\ge Q_n^2(\xi)\int_\eta^{z_j}w\,dt.
$$
[For the proof one considers separately the cases $\eta<\xi$ and
$\eta\ge\xi$.] Using this inequality and the estimate
$\alpha_j=O(\beta)$, we have
$$
{\rm l.h.s.~of~}(3.11)\ge
{Q_n^2(\xi)\int_\eta^{z_j}w\,dt \over
(\min_{F_j}Q_n+\alpha_j)^2\int_{z_{j-1}}^\eta w\,dt}
\ge
\left({\min_{F_j}Q_n+\theta \over \min_{F_j}Q_n+K\beta}\right)^2
(1/\beta-1).
$$
We let $A=(1-c_1+\theta)^2/(1-c_1+\theta/2)^2>1$; then for
$\beta=\theta/2K$ the l.h.s. of (3.11)$\ge\gamma=A(1/\beta-1)$. $\square$

The next lemma will be used in an optimization procedure when proving Lemma 3.6. 
The values of parameters $\beta$ and $\gamma$ will then be taken from Lemma 3.4(ii).

{\Th Lemma 3.5.}{\it For $0<a\le 1/(1+\gamma) < \beta \le 1$ set
$$
F(a;\alpha_1,\alpha_2)=\frac{a}{\alpha_1}+\frac{1-a}{\alpha_2},\quad
\alpha_1,\alpha_2>0.
$$
There exists a constant $M(\beta,\gamma)<1$ such that
$$
\min_{\beta\alpha_1+(1-\beta)\alpha_2=1} F(a;\alpha_1,\alpha_2)\le
  M(\beta,\gamma).
\eqno (3.12)
$$

}
{\Proof Proof.}
The Lagrange function has the form
$$
L(a;\alpha_1,\alpha_2,\lambda)=\frac a{\alpha_1}+\frac{1-a}{\alpha_2}+
\lambda(\beta\alpha_1+(1-\beta)\alpha_2-1).
$$
From $\partial L/\partial\alpha_1=-a/\alpha_1^2+\lambda\beta=0$,
$\partial L/\partial\alpha_2=-(1-a)/\alpha_2^2+\lambda(1-\beta)=0$,
$\beta\alpha_1+(1-\beta)\alpha_2=1$
we have
$$
\alpha_1^*=\left(\frac a{\lambda\beta}\right)^{1/2},\quad
\alpha_2^*=\left(\frac{1-a}{\lambda(1-\beta)}\right)^{1/2},
$$
$$
\lambda=\left((a\beta)^{1/2}+((1-a)(1-\beta))^{1/2}\right)^2.
$$
By algebraic transformations we find the minimum
$$
 \frac a{\alpha_1^*}+\frac{1-a}{\alpha_2^*}=
\lambda^{1/2}\left((a\beta)^{1/2}+((1-a)(1-\beta))^{1/2}\right)=
$$
$$
=1-\left(((1-a)\beta)^{1/2}-(a(1-\beta))^{1/2}\right)^2=
m(\beta,a)\le m\left(\beta,{1 \over 1+\gamma}\right)<1.
$$
Therefore (3.12) holds with
$M(\beta,\gamma)=m(\beta,1/(1+\gamma))$. $\square$

The next lemma is crucial in the proof of Theorem 1.1.
Constants $\varepsilon,r_0$ in Lemmas 3.6 and 3.6$'$ depend only on
$D(u^{-2}),\delta, n$.

{\Th Lemma 3.6.}{\it Let $u^{-2}\in B_\delta$, $n\in\N$. Then there exist constants
$\varepsilon>0$, ${r_0\ge0}$$(=0$ for $\delta=0)$ such that for each
 $r\ge r_0$ one can choose a function $f_r$ on
 $[0,r]$ satisfying the conditions
$$
\|f_r u\|_{L_2(0,r)}<\infty,    \eqno (3.13)
$$
$$
\int_0^r t^k f_r(t)\,dt=0,\quad k=1\ldots n, \eqno (3.14)
$$
$$
\int_0^r f_r(t)\,dt\ge\varepsilon
\|f_r u\|_{L_2(0,r)}\cdot\|u^{-1}\|_{L_2(0,r)}.   \eqno (3.15)
$$

}
{\Proof Proof.} We will look for $f_r$ in the form $u^{-1}g_r$, where
$g_r\in L_2(0,r)$. Conditions (3.14),(3.15) then take the form:
$$
\int_0^r g_r(t)t^k u(t)^{-1}\,dt=0,\quad k=1\ldots n, \eqno (3.16)
$$
$$
\int_0^r g_r(t)u(t)^{-1}\,dt\ge\varepsilon
\|g_r\|_{L_2(0,r)}\cdot\|u^{-1}\|_{L_2(0,r)}.   \eqno (3.17)
$$
Denote by
$\varphi_r=Q_n u^{-1}$, $Q_n(t)=\sum_{k=1}^n b_kt^k$,
the projection in $L_2(0,r)$ of the function $u^{-1}$ on the linear span $E_r$ of the set of functions $\{ t^ku^{-1}: k=1\ldots n\}$.
Let ${g_r=u^{-1}-\varphi_r}=P_n u^{-1}$, $P_n=1-Q_n$.
Condition (3.16) will be clearly satisfied, while
(3.17) will be equivalent to the condition
$$
\exists c<1:\quad\int_0^r \varphi_r(t)u^{-1}(t)\,dt\le c
\|\varphi_r\|_{L_2(0,r)}\cdot\|u^{-1}\|_{L_2(0,r)}.   \eqno (3.18)
$$
Indeed, in terms of geometry of the Hilbert space $L_2(0,r)$
conditions (3.16),(3.17) mean that the angle between the vectors $u^{-1}$ and $g_r$, $
g_r\bot E_r$, is uniformly in $r\ge r_0$ ``small'' (separated from $\pi/2$),
while (3.18) shows that the angle between $u^{-1}$ and $E_r$ is uniformly ``large'' (separated from 0). Clearly these are equivalent statements.

Therefore, we will be proving (3.18).
Condition (3.16) shows that $P_n$ is the $n$-the polynomial of the orthogonal system of polynomials with weight $tu^{-2}$ on the interval $[0,r]$. Since $u^{-2}\in B_\delta$, we have
$tu^{-2}\in B_\delta$ (Lemma 3.1).
By Lemma 3.6 for $w=u^{-2}, v=tu^{-2}$ we can find constants
$\beta\in(0,1)$, $\gamma>1/\beta-1$ such that
for each $r\ge r_0$ there is a set
$A_r\subset[0,r]$ with the properties
$$
\int_{A_r}u^{-2}\,dt=\beta\int_0^r u^{-2}\,dt,  \eqno (3.19)
$$
$$
\int_{[0,r]\backslash A_r}Q_n^2(t)u^{-2}(t)\,dt\ge \gamma
\int_{A_r}Q_n^2(t)u^{-2}(t)\,dt.     \eqno  (3.20)
$$
Let us rewrite (3.20) in the form
$$
\int_{[0,r]\backslash A_r}|\varphi_r(t)|^2\,dt\ge \gamma
\int_{A_r}|\varphi_r(t)|^2\,dt.     \eqno  (3.21)
$$
Consider on $[0,r]$ the function
$$
\alpha(t)=\begin{cases}
     \alpha_1,& t\in A_r, \\
	  \alpha_2, &t\in[0,r]\backslash A_r,\end{cases}
$$
$$
\alpha_1,\alpha_2>0,\quad \beta\alpha_1+(1-\beta)\alpha_2=1.
\eqno (3.22)
$$

By Cauchy-Bunyakovsky-Schwarz inequality and by (3.19)
$$
\left(
\int_0^r \varphi_r(t)u^{-1}(t)\,dt
\right)^2=\left(\int_0^r \frac{\varphi_r(t)}
{\alpha^{1/2}(t)}u^{-1}(t)\alpha^{1/2}(t)\,dt\right)^2\le
$$
$$
\le\int_0^r\frac{|\varphi_r(t)|^2}{\alpha(t)}\,dt
\int_0^ru^{-2}(t)\alpha(t)\,dt =
$$
$$
=\left(\frac1{\alpha_1}\int_{A_r}|\varphi_r(t)|^2\,dt+
\frac1{\alpha_2}\int_{[0,r]\backslash A_r}|\varphi_r(t)|^2\,dt\right)
\int_0^ru^{-2}(t)\,dt.
$$
Let us minimize the last expression over all
$\alpha_1,\alpha_2$ satisfying (3.22).
Applying Lemma 3.5 with $a=\left.\int_{A_r}|\varphi_r(t)|^2\,dt\right/
\int_0^r|\varphi_r(t)|^2\,dt$ [condition $a\le1/(1+\gamma)$ follows
from (3.21)], we have
$$
\min_{\beta\alpha_1+(1-\beta)\alpha_2=1}
\frac1{\alpha_1}\int_{A_r}|\varphi_r(t)|^2\,dt+
\frac1{\alpha_2}\int_{[0,r]\backslash A_r}|\varphi_r(t)|^2\,dt\le
$$
$$
\le M(\beta,\gamma)\int_0^r|\varphi_r(t)|^2\,dt ,\quad
M(\beta,\gamma)<1,
$$
from where (3.18) follows with $c=(M(\beta,\gamma))^{1/2}$. $\square$

The next statement is essentially a reformulation of Lemma 3.6 and easily follows from the remark made after (3.18). Nevertheless we believe that it may be of independent interest.

{\Th Lemma 3.6$'$.}{\it Let $u^{-2}\in B_\delta$, $n\in{\hbox{\bf N}}$, and
$\chi_r=\chi_{(0,r)}$ be the characteristic function of the interval. Then there exist constants
$\varepsilon>0$, $r_0\ge0(=0$ for $\delta=0)$ such that for any $r\ge r_0$
$$
\hbox{{\bf G}}(u^{-1}\chi_r,xu^{-1}\chi_r,\ldots,x^nu^{-1}\chi_r)\ge
$$
$$
\ge\varepsilon\|u^{-1}\chi_r\|_2\cdot\|xu^{-1}\chi_r\|_2\cdots
\|x^nu^{-1}\chi_r\|_2.
$$
}
Here $\hbox{\bf G}$ is the Gram determinant of a system of functions in $L_2$.
The claim, therefore, is that the parallelepiped with edges $u^{-1}\chi_r,\ldots,x^nu^{-1}\chi_r$ is uniformly
in $r$ non-degenerate.

\newpage
\punkt{4. Proofs of results from Section 1}
{\Th Lemma 4.1.}
{\it Let $u^{-2}\in B_\delta$. If $\delta>0$,
assume in addition $a_0v\in L_2(0,r)\, \forall r>0$. Then inequality
{\rm(1.2)} implies that
$$
S_0=\sup_{r>0} \|a_0 v\|_{L_2 (r,\infty)}\cdot
		   \|u^{-1}\|_{L_2 (0,r)} <\infty.
$$

}
{\Proof Proof.} Let us apply Lemma 3.6. For a function $f_r,r\ge r_0$,
satisfying condition (3.13)--(3.15) and extended by zero on
$[r,\infty)$, we will consecutively have
$$
 \| v {\cal A}f_r \|_2 \le C \|  uf_r \|_2 \Rightarrow
\left \| v(x)\sum_{k=0}^na_k(x)\int_0^rt^kf_r(t)\,dt \right\|_{L_2(r,\infty)} \le
  C \|  uf_r \|_{L_2(0,r)}
$$
$$
 \mathop{\Longrightarrow}\limits^{(3.14),(3.15)}
\|a_0v\|_{L_2(r,\infty)}\cdot\|u^{-1}\|_{L_2(0,r)}\le C/\varepsilon,
\quad r\ge r_0.
$$
It remains to note that $r_0=0$ for $\delta=0$, while in the case $\delta>0$ for
$r<r_0$ we have
$$
 \|a_0 v\|_{L_2 (r,\infty)}\cdot\|u^{-1}\|_{L_2 (0,r)} \le
$$
$$
\le (\|a_0 v\|_{L_2 (0,r_0)}+\|a_0 v\|_{L_2 (r_0,\infty)})
 \|u^{-1}\|_{L_2 (0,r_0)} \le
$$
$$
\le \|a_0 v\|_{L_2 (0,r_0)}\cdot\|u^{-1}\|_{L_2 (0,r_0)}+C/\varepsilon.
\quad\square
$$

{\Th Lemma 4.2.} {\it For the inequality 
$$
\left\|v(x)\int_0^xf(t)dt\right\|_2\le c\|uf\|_2
$$ to hold with a constant $c>0$ independent of function $f$,  it is necessary and sufficient that
$$
\sup_{r>0}\|v\|_{L_2 (r,\infty)}\cdot
		   \|u^{-1}\|_{L_2 (0,r)} <\infty.
                   $$
}
This is the known criterion of the weighted Hardy inequality, obtained in [3--5]

{\Proof Proof of Theorem 1.1.} By Lemma 4.2 the condition $S_k<\infty$
is necessary and sufficient for the inequality
$$
\left\|v(x)a_k(x)\int_0^xf(t)dt\right\|_2\le c\|x^{-k}u(x)f(x)\|_2,
$$
which by replacing $\tilde f(x)=x^{-k}f(x)$ becomes
$$
\left\|v(x)a_k(x)\int_0^xt^k\tilde f(t)dt\right\|_2\le c\|
u(x)\tilde f(x)\|_2.
$$
The latter inequality means that
${\cal A}_k:L_{2,u}\to L_{2,v}$. We conclude the equivalence
of Theorems 1.1 and 1.1$'$
and sufficiency of conditions (1.3) of Theorem 1.1.

To prove that conditions (1.3) are necessary we will use the induction on $n$.
For $n=0$ the statement of the theorem follows from Lemma 4.2.
Assume the theorem is proved for $n\le n_0$, and for the operator
$$
  ({\cal A} f)(x)=\int_0^x \left[\sum_{k=0}^{n_0+1} a_k (x)t^k \right]f(t)\,dt.
$$
inequality (1.2) holds. Then Lemma 4.1 implies $S_0<\infty$,
which is equivalent to the inequality
$$
\left\|va_0\int_0^xf\,dt\right\|_2\le c\|uf\|_2.    \eqno(4.1)
$$
From (1.2) and (4.1) we obtain the inequality
$$
\left\|v(x)\int_0^x \left[\sum_{k=1}^{n_0+1} a_k (x)t^k \right]f(t)\,dt\right
\|_2
\le C'\|uf\|_2,
$$
which by substituting $\tilde f(t)=tf(t)$ reduces to the form
$$
\left\|v(x)\int_0^x \left[\sum_{k=0}^{n_0} a_{k+1} (x)t^k \right]
\tilde f(t)\,dt\right\|_2\le C'\|t^{-1}u(t)\tilde f(t)\|_2.
$$
Since $u^{-2}\in B_\delta$, then $(t^{-1}u)^{-2}=t^2u^{-2}\in B_\delta$
(Lemma 3.1), and by the inductive hypothesis we obtain finiteness of the other constants $S_k$ $(k=1\ldots n_0+1)$. $\square$

{\Rem Remark.} For $n=1$ Lemmas 3.6, 4.1, and therefore Theorem 1.1 as well, remain true if
the condition $u^{-2}\in B_\delta$ is replaced in their formulation by the weaker condition (1.4).
In this case Lemma 3.4 used in the proof of Lemma 3.6 is replaced by the following analogously proven result.

{\Th Lemma 3.4$'$.} {\it Let function $w\ge 0$ satisfy with some constants
 $D,\delta$ the condition
$$
\int_0^{2r}w\,dt\le D\int_0^rw\,dt,\quad r\ge\delta\ge0.
$$
Then there exists such constants $\beta\in(0,1)$, $\gamma>1/\beta-1$,
$r_0\ge0(=0$ in the case
$\delta=0)$, that for each $r\ge r_0$
$$
\int_{r^*}^rt^2w(t)\,dt\ge \gamma
\int_0^{r^*}t^2w(t)\,dt,
$$
where $r^*\in(0,r)$ is determined by the condition
$$
\int_0^{r^*}w\,dt=\beta\int_0^r w\,dt.
$$
}

{\Th Lemma 4.3.} {\it Let $\alpha\ge1$. To have the inequality
$\|v\cR^{(\alpha)}f\|_2\le c\|uf\|_2$
it is necessary and sufficient that the two conditions hold:
$$
\sup_{r>0}\|(x-r)^{\alpha-1}v\|_{L_2(r,\infty)}\cdot
		  \|u^{-1}\|_{L_2(0,r)}<\infty ,
$$
$$
\sup_{r>0}\|v\|_{L_2(r,\infty)}\cdot
		  \|(r-x)^{\alpha-1}u^{-1}\|_{L_2(0,r)}<\infty.
$$
}
This is the criterion of boundedness of the Riemann-Liouville operators in weighted spaces, obtained in [6--8].

{\Proof Proof of Assertion 1.2.} Apply Lemma 4.3. $\square$

{\Proof Proof of Assertion 1.3.} By Lemma 4.2 condition (1.6)
is necessary and sufficient for the inequality
$$
\left\|v(x)x^{\alpha-1}\int_0^xf(t)dt\right\|_2\le c\|uf\|_2.
$$
Since $x^{\alpha-1}\ge(x-t)^{\alpha-1}$, $x\ge t\ge 0$, the sufficiency of
(1.6) for (1.5) is clear.

The estimate
$$
 \|x^{\alpha-1}v\|_{L_2(r,\infty)}\cdot
		  \|u^{-1}\|_{L_2(0,r)}\le
$$
$$
\le 2^{\alpha-1}\|(x-r/2)^{\alpha-1}v\|_{L_2(r/2,\infty)}\cdot
		  \|u^{-1}\|_{L_2(0,r)}\le
$$
$$
\le 2^{\alpha-1}D^{1/2}\|(x-r/2)^{\alpha-1}v\|_{L_2(r/2,\infty)}\cdot
		  \|u^{-1}\|_{L_2(0,r/2)}
$$
shows that necessity of (1.6) can be obtained with the help of Lemma 4.3. Such arguments were used in  \cite{Kuf}. $\square$

\punkt{5. Proofs of results from section 2}
{\Proof Proof of Lemma 2.1.}
$$
m<l:\quad(\varphi g)^{(m)}(x)=\sum_{k=0}^mC_m^k\varphi^{(k)}(x)g^{(m-k)}(x)=
$$
$$
=\sum_{k=0}^mC_m^k\varphi^{(k)}(x)
\int_0^x\frac{(x-t)^{l-m+k-1}}{(l-m+k-1)!}g^{(l)}(t)\,dt=
$$
$$
=\frac1{(l-1)!}\int_0^x\left[\sum_{k=0}^mC_m^k\varphi^{(k)}(x)
\frac{(l-1)!}{(l-m+k-1)!}(x-t)^{l-m+k-1}\right]g^{(l)}(t)\,dt=
$$
$$
=\frac1{(l-1)!}\int_0^x\frac{d^m}{dx^m}(\varphi(x)(x-t)^{l-1})
g^{(l)}(t)\,dt.
\eqno (5.1)
$$
$$
(\varphi g)^{(l)}(x)=\varphi(x)g^{(l)}(x)+
\sum_{k=1}^lC_l^k\varphi^{(k)}(x)g^{(l-k)}(x)=
$$
$$
=\varphi(x)g^{(l)}(x)+\sum_{k=1}^{l}\varphi^{(k)}(x)
\int_0^x{(x-t)^{k-1} \over (k-1)!}g^{(l)}(t)\,dt=
$$
$$
=\varphi(x)g^{(l)}(x)+\frac1{(l-1)!}
\int_0^x\left[\sum_{k=1}^{l}C_{l}^k\varphi^{(k)}(x)
{(l-1)! \over (k-1)!}(x-t)^{k-1}\right]g^{(l)}(t)\,dt=
$$
$$
=\varphi(x)g^{(l)}(x)+\frac1{(l-1)!}
\int_0^x{d^l \over dx^l}(\varphi(x)(x-t)^{l-1})g^{(l)}(t)\,dt.
\eqno (5.2)
$$
To get (2.2) and (2.3) it remains to expand $(x-t)^{l-1}$
in (5.1) and (5.2) by the binomial formula. $\square$

{\Th Lemma 5.1.}{\it For any set $h_1,...,h_l$ of functions integrable
on $[a,b]$, one can find a function $\sigma$ with $|\sigma(x)|=1$ on $[a,b]$
such that
$$
\int_a^bh_k(x)\sigma(x)\,dx=0,\quad k=1\ldots l.
$$

}
{\Proof Proof} of this statement can be found in the book
\cite[p.267]{KS}. $\square$

{\Proof Proof~of lemma~2.2.} Define the norm $\|\varphi\|_M$ in the space of multipliers
as the norm of the corresponding operator multiplying by $\varphi$,
acting from $W_{2,u}^{(l)}$ into $W_{2,v}^{(l)}$.

From the conditions
$u^{-1},v^{-1}\in L_2(0,r)\,\forall r>0$ one can deduce completeness of the considered
$W$--spaces. In this case from Banach's closed graph theorem it is easy to get that $\|\varphi\|_M<\infty$. Therefore it suffices to prove the inequality
$$
\|\varphi v u^{-1}\|_{L_\infty(\RP)}\le c\|\varphi\|_M. \eqno(5.3)
$$
Let $g\in W_{2,u}^{(l)}$ and satisfies (2.1). By Lemma 2.1
$$
\|\varphi g\|_{W_{2,v}^{(l)}}\ge\|(\varphi g)^{(l)}v\|_2\ge
$$
$$
\ge\|\varphi g^{(l)}v\|_2-c\sum_{k=0}^{l-1}
\left\|(\varphi x^k)^{(l)}v\int_0^xt^{l-k-1}g^{(l)}(t)\,dt\right\|_2.
\eqno (5.4)
$$
Furthermore, consider the set $A_\alpha=\{x:|\varphi v u^{-1}(x)|\ge\alpha\}$.
Let $\mbox{mes } A_\alpha>0$. Then $\forall \varepsilon>0$ there exists an interval
$\Delta_\varepsilon$, $|\Delta_\varepsilon|=\varepsilon$, such that
$\mbox{mes } \Delta_\varepsilon\cap A_\alpha>0$. By Lemma 5.1 there exists a function $g_l$ such that
$$
|g_l(t)|=\begin{cases}
0,  &t\notin \Delta_\varepsilon\cap A_\alpha,\\
u^{-2}(t), &t\in \Delta_\varepsilon\cap A_\alpha,
\end{cases}
$$
and
$$
\int_{\Delta_\varepsilon}t^{l-k-1}g_l(t)\,dt=0,\quad k=0\ldots l-1.
$$
We put
$$
g(x)=\int_0^x\frac{(x-t)^{l-1}}{(l-1)!}g_l(t)\,dt.
$$
Then $g^{(l)}=g_l$, (2.1) is satisfied, while (5.4) gives:
$$
\|\varphi g\|_{W_{2,v}^{(l)}}\ge\alpha\|g^{(l)}u\|_
{L_2(\Delta_\varepsilon)}-
$$
$$
-c\sum_{k=0}^{l-1}\|(\varphi x^k)^{(l)}v\|_{L_2(\Delta_\varepsilon)}\cdot
\|g^{(l)}u\|_{L_2(\Delta_\varepsilon)}\cdot
\|x^{l-k-1}u^{-1}\|_{L_2(\Delta_\varepsilon)}=
\eqno (5.5)
$$
$$
=\left(\alpha-c\sum_{k=0}^{l-1}\|(\varphi x^k)^{(l)}v\|_
{L_2(\Delta_\varepsilon)}\cdot
\|x^{l-k-1}u^{-1}\|_{L_2(\Delta_\varepsilon)}\right)
\|g^{(l)}u\|_{L_2(\Delta_\varepsilon)}.
$$
Note that $\|g\|_{L_2(0,1)}\le\|g^{(l)}u\|_2\cdot\|u^{-1}\|_{L_2(0,1)}$, therefore
$\|g\|_{W_{2,u}^{(l)}}\le c\|g^{(l)}u\|_2$.
Since $\|g^{(l)}u\|_2>0$, we get $\|\varphi\|_M\ge\alpha/c$ from (5.5) by tending $\varepsilon$ to zero. Tending now $\alpha$ to
$\mbox{ess sup}_{x>0}|\varphi v u^{-1}(x)|-0$, we get (5.3). $\square$

{\Proof Proof of Theorems 2.3 and 2.4.} Consider the case $m=l$. The case
$m<l$
is considered analogously

{\it Necessity.} Let $\varphi\in M(u,l;v,l)$. Then
(2.5) and (2.8) follow from Lemma 2.2, while (2.4) and (2.6) follow from
$x^k\in W_{2,u}^{(l)}$, $k=0\ldots l-1$. It remains to show the necessity of
condition (2.7) in Theorem 2.4. From the inequality
$$
\|(\varphi g)^{(l)}v\|_2\le c\|g^{(l)}u\|_2
$$
considered on functions satisfying (2.1), Lemma 2.1 and (2.8)
we get the inequality
$$
\left\|\sum_{k=0}^{l-1}C_{l-1}^k(\varphi x^k)^{(l)}v
\int_0^x(-t)^{l-k-1}g^{(l)}(t)\,dt\right\|_2\le c
\|g^{(l)}u\|_2.
$$
Now (2.7) follows from Theorem 1.1.

{\it Sufficiency. } Let us show that if
(2.6)--(2.8) hold and $f\in W_{2,u}^{(l)}$, then $\varphi f\in W_{2,v}^{(l)}$.
This will show sufficiency of conditions of both Theorem 2.4 and Theorem 2.3, since in the case $(1+x^{l-1})u^{-1}\in L_2$ conditions (2.6),(2.8)
imply (2.7).

Representing $f$ in the form
$$
f(x)=\sum_{k=0}^{l-1}\frac{x^k}{k!}f^{(k)}(0)+g(x),
$$
where $g$ clearly satisfies (2.1), we will have
$$
\|\varphi f\|_{W_{2,v}^{(l)}}\le
c\sum_{k=0}^{l-1}\|\varphi x^k\|_{W_{2,v}^{(l)}}+
\|\varphi g\|_{W_{2,v}^{(l)}}=
$$
$$
=c\sum_{k=0}^{l-1}\left(\|\varphi x^k\|_{L_2(0,1)}+
\|(\varphi x^k)^{(l)}v\|_2\right)+\|\varphi g\|_{L_2(0,1)}+
\|(\varphi g)^{(l)}v\|_2.
$$
With the help of Lemmas 2.1 and 4.2, (2.7) and (2.8) imply the inequality
$$
\|(\varphi g)^{(l)}v\|_2\le c\|g^{(l)}u\|_2,
$$
providing an estimate for the last term. Taking into account (2.6) it remains to show that
$$
\|\varphi x^k\|_{L_2(0,1)},\|\varphi g\|_{L_2(0,1)}<\infty.
$$
But this is obvious, since $\|\varphi^{(l)}v\|_2,\|g^{(l)}u\|_2<\infty$
imply the continuity of $\varphi$ and $g$ on $\RP$. $\square$

\end{document}